\documentclass[11pt]{amsart}
\usepackage{amscd}
\usepackage[arrow,matrix]{xy}

\theoremstyle{plain}

\newtheorem{thm}{Theorem}[section]
\newtheorem{theorem}[thm]{Theorem}

\newtheorem{lemma}[thm]{Lemma}

\newtheorem{prop}[thm]{Proposition}
\newtheorem{proposition}[thm]{Proposition}
\newtheorem{conjecture}[thm]{Conjecture}
\newtheorem*{conjecture*}{Conjecture}

\newtheorem*{question*}{Question}

\theoremstyle{definition}
\newtheorem{defn}[thm]{Definition}

\newtheorem*{rem*}{Remark}
\newtheorem{remark}[thm]{Remark}
\newtheorem*{remark*}{Remark}

\newtheorem*{remarks*}{Remarks}

\newtheorem*{example*}{Example}
\newtheorem{example}[thm]{Example}
\newtheorem*{examples*}{Examples}

\newtheorem*{notation*}{Notation}

\newtheorem*{bibliographical-note}{Bibliographical note}
\newcommand{\zeroindent}{\parindent0cm \parskip1ex}
\newcommand{\acknowledgements}{{\em Acknowledgements. }}

{\begin{list}{(\alph{enumi}) }{\usecounter{enumi}

\leftmargin0cm \labelsep0cm \rightmargin0cm 
\setlength{\labelwidth}{\fill}}}
{\end{list}}

\newenvironment{normallist}%
{

\begin{enumerate}}{\end{enumerate}}

{

\begin{enumerate}}{\end{enumerate}}

{

\begin{enumerate} \parsep0cm 
\itemsep0cm \parskip0cm}{\end{enumerate}}

\newenvironment{primelist}%
{

\begin{enumerate}}{\end{enumerate}}

{

\begin{enumerate}}{\end{enumerate}}

\newcommand{\R}{\mathbb{R}}
\newcommand{\Z}{\mathbb{Z}}
\newcommand{\Q}{\mathbb{Q}}
\newcommand{\C}{\mathbb{C}}
\newcommand{\N}{\mathbb{N}}

\newcommand{\id}{\mathrm{id}}
\newcommand{\im}{\mathrm{im}}

\newcommand{\iso}{\cong}           %isomorphism sign
\newcommand{\htp}{\simeq}            %homotopy sign
\newcommand{\smooth}{C^\infty}
\newcommand{\CP}[1]{\C {\mathrm P}^{#1}}

\newcommand{\mo}{(M,\omega)}
\renewcommand{\o}{\omega}
\renewcommand{\O}{\Omega}
\newcommand{\Aut}{\mathrm{Aut}}

\newcommand{\leftsc}{\langle}
\newcommand{\rightsc}{\rangle}
\newcommand{\Diff}{\mathrm{Diff}}

\theoremstyle{plain}
\newtheorem*{assertion*}{Assertion}
\zeroindent
\numberwithin{equation}{section}

\newcommand{\rank}{\mathrm{rank }}
\newcommand{\field}{\mathbb{F}}
\newcommand{\Def}{\mathrm{Def}}
\newcommand{\tR}{\widetilde{R}}
\newcommand{\tu}{\tilde{u}}
\newcommand{\tv}{\tilde{v}}
\renewcommand{\P}{\mathbb{P}}

\newcommand{\qpa}{\ast_A}
\newcommand{\tqpa}{\tilde{\ast}_A}
\newcommand{\tpsi}{\widetilde{\psi}}
\newcommand{\txi}{\tilde{\xi}}
\newcommand{\point}{point}

\newcommand{\moduli}{\mathcal{M}}
\newcommand{\J}{\mathbf{J}}

\newcommand{\vdim}{\mathrm{v. dim.}}
\newcommand{\coker}{\mathrm{coker} \;}

\title[Automorphisms of $\CP{m} \times \CP{n}$]{On the group of 
symplectic\\ automorphisms of $\CP{m} \times \CP{n}$}
\author{Paul Seidel}
\date{19/3/1998}
\thanks{Supported by NSF grant DMS 9304580.}
\address{Institute for Advanced Study, Princeton\\}
\email{pseidel@math.ias.edu}
\begin{document}
\maketitle

%------------------\input{1}
\section{Introduction \label{sec:one}}

Let $\o_n$ be the standard symplectic form on $\CP{n}$, normalized
in such a way that $[\o_n]$ is Poincar{\'e} dual to a hyperplane.
We denote the product $(\CP{m},\o_m) \times (\CP{n},\o_n)$, for
$m,n \geq 1$, by $(P_{mn},\eta_{mn})$. Let $\Diff(P_{mn})$ be the
group of diffeomorphisms with the $\smooth$-topology,
$\Aut(P_{mn},\eta_{mn})$ the subgroup of symplectic automorphisms, and
\[
\beta_k: \pi_k(\Aut(P_{mn},\eta_{mn})) \longrightarrow
\pi_k(\Diff(P_{mn}))
\]
the homomorphisms induced by inclusion.

\begin{theorem} \label{th:main} Let $k$ be odd and 
$\leq \max\{2m-1,2n-1\}$. Then $\beta_k$ is not surjective.
In fact
\begin{multline*}
\rank(\coker \beta_k) \geq 
b_{2m+1-k}(P_{mn}) - b_{2m+1-k}(\CP{m}) + \\
+ b_{2n+1-k}(P_{mn}) - b_{2n+1-k}(\CP{n}) > 0,
\end{multline*}
where $b_*$ are the Betti numbers.
\end{theorem}

For $m = n = 1$ the Theorem says that
\begin{equation} \label{eq:mn-one}
\rank(\coker \beta_1) \geq 2.
\end{equation}
This can be derived from a result of Gromov \cite[0.3.C]{gromov85} 
which says that the group of K{\"a}hler isometries of $P_{11}$ 
(with respect to the standard metric) is a deformation retract of
$\Aut(P_{11},\eta_{11})$. Since the isometry group is an extension
of $\Z/2$ by $SO(3) \times SO(3)$, it follows that 
$\pi_1(\Aut(P_{11},\eta_{11})) \iso \Z/2 \oplus \Z/2$. The topology
of $\Diff(P_{11})$ is unknown, but by looking at its image in
the space of continuous self-maps of $P_{11}$ one can show that
$\pi_1(\Diff(P_{11}))$ is a group of rank $\geq 2$, which implies
\eqref{eq:mn-one}.

Gromov's theorem on $\Aut(P_{11},\eta_1)$ and its cousing for
$(\CP{2},\o_2)$ were the first non-trivial results about the
topology of symplectic automorphism groups. More recently Abreu
\cite{abreu97} and Abreu-McDuff \cite{abreu-mcduff98} have
studied the group of automorphisms of the symplectic structures
$\eta_{11}^{(\lambda)} = \lambda(\o_1 \times 1) + 1 \times \o_1$
on $P_{11}$ for $\lambda \neq 1$, using an approach suggested by
Gromov \cite[2.4.$C_2$]{gromov85}. The automorphism groups of
blowups of $\CP{2}$ are closely related to the spaces of symplectic
embeddings of balls, which have been studied by Biran \cite{biran94}
and McDuff \cite{mcduff96}. These results rely on specific properties
of rational or ruled symplectic four-manifolds, and little is known
outside this class (although see 
\cite{seidel97} for some information about $\pi_0$ of the automorphism
groups of symplectic four-manifolds). It seems that 
Theorem \ref{th:main} provides the first examples of symplectic manifolds of
dimension $>4$ for which the map $\pi_1(\Aut) \longrightarrow \pi_1(\Diff)$
is known not to be surjective. The fact that no such examples were known
was pointed out to the author by F. Lalonde. Lalonde also suggested that 
it may be possible to construct such examples by exploiting the `rigidity 
theorem' of \cite{lalonde-mcduff-polterovich97} instead of the methods used here 
(any such examples would have positive first Betti number).

By definition, the symplectic automorphism group is the structure group of
symplectic fibre bundles. The study of pseudo-holomorphic curves in each
fibre of such a fibre bundle $E$ yields `parametrized Gromov-Witten invariants'
which are multilinear maps on $H^*(E)$. These
invariants have been considered by L{\^e} \cite{le97} and others. They 
satisfy (at least in principle) axioms similar to those of Kontsevich and 
Manin \cite{kontsevich-manin94}. In this paper we argue that in certain
cases the existence of such a system of invariants restricts the 
possibilities for what $H^*(E)$ can be. 
To simplify the technical issues, we do not consider all Gromov-Witten
invariants but just a single particularly simple one.

Although our argument relies on Gromov-Witten invariants, it does not
suppose any knowledge of what the actual value of the invariants is.
In this respect it resembles the recent `rigidity theorem' of 
Lalonde-McDuff-Polterovich \cite{lalonde-mcduff-polterovich97}. Note that 
the proof of the `rigidity theorem' involves pseudo-holomorphic curves of 
a different kind, namely, pseudo-holomorphic sections of a fibre bundle 
whose base is a Riemann surface (in that case, $S^2$). It is possible that 
a combination of the two approaches would yield more information about
symplectic automorphism groups.

This paper is structured as follows: the next section reviews some
basic facts about the cohomology rings of fibrations. In section 
\ref{sec:three} we apply these considerations to the case where
the fibre is $P_{mn}$. Up to this point the argument is purely 
topological. The relevant Gromov-Witten invariant is introduced in 
sections \ref{sec:four} and \ref{sec:five}. Section \ref{sec:six} 
contains the main computation. In the final section we derive
Theorem \ref{th:main} and discuss a related result and some
possible further developments.

\acknowledgements The author is indebted to Fran{\c c}ois Lalonde,
Dusa McDuff, and Richard Thomas for helpful discussions.

%------------------\input{2}

\section{The cohomology rings of fibrations over spheres}

Fix a field $\field$. By a graded $\field$-algebra 
$R = \bigoplus_{i \in \Z} R^i$ we
mean one which is finite-dimensional, commutative (in the
graded sense), and has a unit. All homomorphisms between
graded algebras preserve units.

\begin{defn} Let $R$ be a graded 
$\field$-algebra. A deformation of $R$ of dimension $d>0$ consists of
\begin{normallist}
\item a graded $\field$-algebra $\tR$,
\item \label{item:flat}
an element $t \in \tR^d$ with $t^2 = 0$, such that the sequence
\[
0 \longrightarrow \tR/t \tR \stackrel{\cdot t}{\longrightarrow}
\tR \longrightarrow \tR/t \tR \longrightarrow 0
\]
is exact, and
\item a homomorphism of graded algebras
$j: \tR \longrightarrow R$ which is surjective with 
kernel $t \tR$.
\end{normallist}
\end{defn}

To be precise, these objects should be called `first order
infinitesimal deformations'; we use the shorter name for
brevity's sake. The exactness of the sequence 
$\tR/t \tR \longrightarrow \tR \longrightarrow \tR/t \tR$ 
is equivalent to the flatness of $\tR$ as a module over 
$\field[\epsilon]/\epsilon^2$, where $\epsilon$ acts by
multiplication with $t$.

Given two graded algebras $R_1, R_2$ and a homomorphism
$f: R_1 \longrightarrow R_2$, one defines a {\em morphism
over $f$} from a $d$-dimensional deformation
$(\tR_1,t_1,j_1)$ of $R_1$ to a $d$-dimensional deformation 
$(\tR_2,t_2,j_2)$ of $R_2$ to be a homomorphism of graded algebras
$\tilde{f}: \tR_1 \longrightarrow \tR_2$ such that 
$\tilde{f}(t_1) = t_2$ and $j_2 \tilde{f} = f j_1$. 
In the special case $R_1 = R_2 = R$ and $f = \id$,
$\tilde{f}$ is called a morphism of deformations of $R$.
All morphisms of deformations of $R$ are isomorphisms.

The set $\Def_d(R)$ of isomorphism classes of 
$d$-dimensional deformations of $R$ carries a natural
structure of an $\field$-vector space. The sum of two $d$-dimensional
deformations $(\tR_i,t_i,j_i)$, $i = 1,2$, is the deformation
$(\tR,t,j)$ defined as follows: let $\Delta \subset \tR_1 \oplus
\tR_2$ be the subalgebra of those $(x_1,x_2)$ such that
$j_1(x_1) = j_2(x_2)$. Then $\tR = \Delta / (t_1,-t_2)\Delta$, 
$t = [t_1,0] = [0,t_2] \in \tR$, and $j: \tR \longrightarrow R$ is
given by $j([x_1,x_2]) = j_1(x_1) = j_2(x_2)$.

\begin{remark} \label{th:explicit} There is a 
second and more explicit definition of $\Def_d(R)$. Let
$Z_d(R)$ be the space of $\field$-bilinear maps
$\psi: R \times R \longrightarrow R$ of degree $-d$
which are graded symmetric and satisfy
\[
(-1)^{d \deg(x)} x\psi(y,z) - \psi(xy,z) + \psi(x,yz) - \psi(x,y)z = 0,
\]
and $B_d(R) \subset Z_d(R)$ the subset of maps of the form
\[
\psi(x,y) = (-1)^{d \deg(x)}x\xi(y) - \xi(xy) + \xi(x)y
\]
where $\xi: R \longrightarrow R$ is some linear map of
degree $-d$. Then $\Def_d(R) = Z_d(R)/B_d(R)$. The
equivalence of this definition with the previous one
is proved by choosing, for a deformation $(\tR,t,j)$ of
$R$, an isomorphism of $\tR$ with $R[\epsilon]/\epsilon^2$ as an 
$\field[\epsilon]/\epsilon^2$-module. Then the product
on $\tR$ has the form
\[
(x_0 + \epsilon x_1) (y_0 + \epsilon y_1) = x_0y_0 + 
\epsilon((-1)^{d \deg(x_0)}x_0y_1 + x_1y_0 + \psi(x_0,y_0))
\]
for some $\psi \in Z_d(R)$. The equivalence class of $\psi$ in 
$Z_d(R)/B_d(R)$ is independent of the choice of isomorphism. 
This alternative description shows that if $R$ is concentrated in 
even dimensions then $\Def_d(R) = 0$ for all odd $d$.
\end{remark}

In the deformation theory of algebras (see  
\cite{gerstenhaber-schack90} for a survey)
$\Def_*(R)$ is called the second Harrison 
cohomology group of the graded algebra $R$.

Let $M$ be a compact manifold. A smooth fibre bundle over 
a sphere $S^d$ with fibre $M$ consists of a manifold $E$, 
a submersion $\pi: E \longrightarrow S^d$,
and a diffeomorphism $i: M \longrightarrow E_{b_0}$ for
some fixed $b_0 \in S^d$. The cohomology of $E$ (unless 
otherwise stated, all cohomology groups have $\field$-coefficients)
sits in a Wang sequence
\[
\dots \stackrel{\delta}{\longrightarrow}
H^*(M) \stackrel{i_!}{\longrightarrow} H^{*+d}(E)
\stackrel{i^*}{\longrightarrow} H^{*+d}(M) 
\stackrel{\delta}{\longrightarrow} H^{*+1}(M)
\longrightarrow \dots
\]
where $i_!$ is the cohomology transfer or pushforward.
Let $\epsilon$ be the standard generator of $H^d(S^d)$,
and $t = \pi^*(\epsilon) \in H^d(E)$.

\begin{lemma} \label{th:deformation}
If $i^*$ is onto then $(H^*(E),t,i^*)$
is a $d$-dimensional deformation of $H^*(M)$.
\end{lemma}

\proof Because $t$ is Poincar{\'e} dual to the fibre of $E$,
$i_!(i^*(x)) = x t$ for all $x \in H^*(E)$. Hence $\ker(i^*) = \im(i_!) =
\im(i_!i^*) = t H^*(E)$. It remains to show that any $x$ with $tx = 0$
lies in $tH^*(E)$. Because $i_!$ is injective, $tx = i_!i^*(x) = 0$
implies that $i^*(x) = 0$, and we have already seen that $\ker(i^*) =
t H^*(E)$. \qed

The isomorphism classes of smooth fibre bundles over $S^d$ with fibre $M$ 
form a group under the operation of fibre connected sum. This
group is isomorphic to $\pi_{d-1}(\Diff(M))$; the isomorphism
is given by a clutching construction which associates to
a map $\phi: (S^{d-1},b_0) \longrightarrow (\Diff(M),\id)$
a fibre bundle $(E_\phi,\pi_\phi,i_\phi)$. On the level
of cohomology $\phi$ and $E_\phi$ are related in the
following way: since $\phi$ can be written as a map
$M \times S^{d-1} \longrightarrow M$, it induces a
homomorphism $\delta_\phi: H^{*+d-1}(M) \longrightarrow H^*(M)$,
and this homomorphism is the connecting map in the Wang sequence 
for $E_\phi$. In particular, the subgroup $G_{d-1} \subset
\pi_{d-1}(\Diff(M))$ of classes $[\phi]$ such that $\delta_\phi = 0$
corresponds to the isomorphism classes of fibre bundles $(E,\pi,i)$ 
for which $i^*$ is onto. Lemma \ref{th:deformation} associates to any 
such fibre bundle an element of $\Def_d(H^*(M))$. Since the sum of
deformations imitates the behaviour of cohomology under fibre connected sum,
this defines a homomorphism
\[
\alpha_{d-1}: G_{d-1} \longrightarrow \Def_d(H^*(M)).
\]
Although we have used smooth fibre bundles throughout, the
construction does not really depend on the smooth structure
of $M$. In fact the group of diffeomorphisms can be replaced 
by the bigger semigroup of homotopy equivalences
$M \longrightarrow M$. To take into account the smooth
structure of $M$ would mean to consider deformations
with certain distinguished elements (the Pontryagin
classes). For $\field = \Q$ these questions
can be treated in a much more satisfactory way in the 
framework of rational homotopy theory; see 
\cite[pp. 313--314, 322--326]{sullivan77}.

\begin{example} \label{ex:cpn}
Let $R = \field[u]/u^{n+1}$ where $u$ has dimension $2$. 
Since $R$ is concentrated in even dimensions,
$\Def_d(R) = 0$ for all odd $d$. Now fix an even $d$.
To any element $a \in R^{2n+2-d}$ one can
associate a $d$-dimensional deformation $(\tR_a,t,j)$
of $R$. This deformation is defined in the following way:
$a = \alpha u^{n+1-d/2}$ for some $\alpha \in \field$. Then
\[
\tR_a = \field[\tu,t]/(\tu^{n+1} + \alpha t \tu^{n+1-d/2}, t^2),
\]
and $j: \tR_a \longrightarrow R$ is the algebra homomorphism 
with $j(\tu) = u$ and $j(t) = 0$. The sum of $\tR_a$ and $\tR_b$
is isomorphic to $\tR_{a+b}$.

Every $d$-dimensional deformation $(\tR,t,j)$ of $R$ is 
isomorphic to $\tR_a$ for some $a$. This is proved as follows:
choose a $\tu \in \tR^2$ with $j(\tu) = u$. 
Since $u^{n+1} = 0$, $\tu^{n+1} = -t\tilde{a}$ for 
some $\tilde{a}$. At this point it is easy to construct a morphism of 
deformations of $R$ from $\tR_a$, where $a = j(\tilde{a})$, to $\tR$; 
recall that any such morphism is an isomorphism.

If $d \neq 2$ then $\tR_a$ is isomorphic to $\tR_b$
iff $a = b$. This can be proved either directly or
by observing that for an arbitrary deformation, the
element $\tu$ is unique and hence $a = j(\tilde{a})$ is
an isomorphism invariant of the deformation. For $d = 2$
$\tR_a$ and $\tR_b$ are isomorphic iff $b-a$ is a multiple
of $(n+1)u^n$. Therefore
\[
\Def_d(R) \iso \begin{cases}
H^{2n+2-d}(R) & d \neq 2,\\
H^{2n}(R)/(n+1)H^{2n}(R) & d = 2.
\end{cases}
\]
We will now interpret this result geometrically. $R = H^*(\CP{n})$. 
Because the cohomology is concentrated in even dimensions, any smooth fibre 
bundle $(E,\pi,i)$ with fibre $\CP{n}$ over $S^d$ satisfies the condition 
that $i^*$ is onto if $d$ is even. In other words $G_{d-1} = 
\pi_{d-1}(\Diff(\CP{n}))$ for all even $d$. Hence one has homomorphisms
\[
\alpha_{d-1}: \pi_{d-1}(\Diff(\CP{n})) \longrightarrow 
\Def_{d}(R) \iso \begin{cases}
\field & d = 4, \dots, 2n+2,\\
\field/(n+1)\field & d = 2.
\end{cases}
\]
If we restrict them to the subgroup $PU(n+1) \subset \Diff(\CP{n})$, these
homomorphisms are the Chern classes $c_2, \dots, c_{n+1}$ and the first 
Chern class mod $n+1$. This follows from the standard formula for the 
cohomology ring of the projective bundle associated to a vector bundle.
\end{example}

The argument above is an instance of a general way of computing the 
deformation spaces of an algebra given by generators and relations.
We will use this method again in the next section.

%------------------\input{3}

\section{Deformations of $H^*(P_{mn})$ \label{sec:three}}

Let $R = H^*(P_{mn}) = \field[u,v]/(u^{m+1},v^{n+1})$. Fix an even $d$. 
One can associate to any pair $(a,b) \in R^{2m+2-d} \oplus R^{2n+2-d}$ 
a $d$-dimensional deformation of $R$ in the following way: write 
$
a = \sum_{i} a_i u^{i-d/2} v^{m+1-i}
$
and
$
b = \sum_{i} b_i u^{i-d/2} v^{n+1-i}
$
with $a_i,b_i \in \field$. Let $\tR_{a,b}$ be the algebra with generators 
$\tu,\tv,t$ of degrees $2,2,d$ and relations $t^2 = 0$,
\[
\tu^{m+1} + \sum_i a_i t \tu^{i-d/2} \tv^{m+1-i} = 0, \quad
\tv^{n+1} + \sum_i b_i t \tu^{i-d/2} \tv^{n+1-i} = 0.
\]
$\tR_{a,b}$, $t$, and the homomorphism $j: \tR_{a,b} \longrightarrow R$ with
$j(\tu) = u$, $j(\tv) = v$, and $j(t) = 0$, define a deformation of $R$.

Using the same idea as in Example \ref{ex:cpn}, one can check that for
$d \neq 2$ these examples form a complete list, with no repetitions, of the
$d$-dimensional deformations of $R$. For $d = 2$ there are again some
isomorphisms between the $\tR_{a,b}$ for different $(a,b)$. In this way
one obtains isomorphisms
\begin{equation} \label{eq:isomorphism}
\Def_d(R) \iso \begin{cases} 
R^{2m+2-d} \oplus R^{2n+2-d} & d \neq 2,\\
R^{2m}/\field(m+1)u^m \oplus R^{2n}/\field(n+1)v^n & d = 2.
\end{cases}
\end{equation}
As in the case of $\CP{n}$, the homomorphisms $\alpha_{d-1}$
are defined on all of $\pi_{d-1}(\Diff(P_{mn}))$ if $d$ is even.
For the rest of this section let $\field = \Q$.

\begin{prop} \label{th:surjective} 
For any even $d$ the homomorphism
\[
\alpha_{d-1} \otimes \id_{\Q}: \pi_{d-1}(\Diff(P_{mn})) \otimes \Q
\longrightarrow \Def_d(R)
\]
is surjective.
\end{prop}

Clearly the proof consists in finding sufficiently many examples of
smooth fibre bundles. We will use the following construction:
let $\xi \longrightarrow \CP{m} \times S^d$ be a complex vector
bundle of rank $n+1$, such that $\xi|\CP{m} \times \{b_0\}$ is
trivial, and $\P(\xi)$ the associated bundle of projective spaces.
The map $\P(\xi) \longrightarrow \CP{m} \times S^d \longrightarrow S^d$
makes $\P(\xi)$ into a smooth fibre bundle with fibre $\CP{m} \times \CP{n}$.
$H^*(\P(\xi))$ is a module over $H^*(\CP{m} \times S^d)$ with one 
two-dimensional generator $v$ and one relation
\[
v^{n+1} + \sum_{i=0}^n c_i(\xi) v^{n+1-i} = 0.
\]
Since $\xi|\CP{m} \times \{b_0\}$ is trivial, the Chern classes can be
written as $c_i(\xi) = \gamma_i (u^{i-d/2} \otimes \epsilon)$, where 
$\gamma_i \in \Z$, $u$ is the generator of $H^2(\CP{m})$, and $\epsilon$
is the generator of $H^2(S^d)$. It follows that the deformation of
$R$ determined by $\P(\xi)$ is isomorphic to $\tR_{0,b}$ with
$b = \sum_i \gamma_i u^{i-d/2} v^{n+1-i}$.
The Chern classes which can be nonzero are those with $d/2 \leq i
\leq \nu = \min\{n+1,m+d/2\}$. Now assume that for every such $i$
there is a vector bundle $\xi_i$ such that $c_i(\xi_i) \neq 0$ and
$c_j(\xi_i) = 0$ for all $j>i$. Since
$R^{2n+2-d} = \bigoplus_{i = d/2}^{\nu} \Q u^{i-d/2} v^{n+1-i}$,
it would follow that $0 \oplus R^{2n+2-d} \subset \Def_d(R)$ lies
in the image of $\alpha_{d-1} \otimes \id_{\Q}$. By exchanging $m$
and $n$, one would obtain the same result for the complementary
subspace $R^{2m+2-d} \oplus 0$, and this would prove Proposition
\ref{th:surjective}. Hence it remains to construct the $\xi_i$. In order
to fulfil $c_j(\xi_i) = 0$ for $j>i$, we will take $\xi_i$ to be the sum 
of a vector bundle $\eta_i$ of rank $i$ and a trivial bundle. All we
need to prove is that

\begin{lemma} For any integer $i$ with $d/2 \leq i \leq m+d/2$ there is a 
complex vector bundle $\eta_i \longrightarrow \CP{m} \times S^d$ of rank $i$
such that $\eta_i|\CP{m} \times \{b_0\}$ is trivial and $c_i(\eta_i) \neq 0$.
\end{lemma}

\proof Vector bundles $\eta \longrightarrow \CP{k} \times S^d$ of rank $i$ 
equipped with a trivialization of $\eta|\CP{k} \times \{b_0\}$ are
classified by the homotopy classes of (based) maps
\begin{equation} \label{eq:maps}
\CP{k} \longrightarrow \Omega^{d-1} U(i) \htp \Omega^d BU(i).
\end{equation}
Since $\Omega^{d-1}U(i)$ is a topological group, two such maps can be
multiplied. This multiplication associates to two vector bundles a 
third one, whose $i$-th Chern class is the sum of the $i$-th Chern classes 
of the original ones.

The obstruction to extending a map \eqref{eq:maps} from $\CP{k}$ to 
$\CP{k+1}$ is an element of the finite group 
$\pi_{2k+1}(\Omega^{d-1} U(i)) = \pi_{2k+d}(U(i))$. Hence the extension 
can always be carried out after replacing the original map by a positive
multiple. Start with a map $\CP{i-d/2} \longrightarrow \Omega^{d-1} U(i)$
which collapses everything except the top-dimensional cell to a point
and represents a nontrivial element of $\pi_{2i-d}(\Omega^{d-1}U(i)) 
\otimes \Q \iso \pi_{2i-1}(U(i)) \otimes \Q \iso \Q$. This map corresponds
to a vector bundle over $\CP{i-d/2} \times S^d$ with nonzero $i$-th
Chern class. After passing to a suitable multiple if necessary, one
can extend the map to $\CP{m}$. \qed

%------------------\input{4}

\section{Symplectic fibre bundles \label{sec:four}}

In this section and the next one, $\mo$ is a compact symplectic manifold, 
$R = H^*(M)$, and $A \in H_2^s(M;\Z) = \im(\pi_2(M) \longrightarrow H_2(M;\Z))$
is a spherical homology class such that $\o(A)$ is positive and 
generates the period group $\o(H_2^s(M;\Z)) \subset \R$. 
In other words, we assume that there is no $A' \in H_2^s(M;\Z)$ 
with $0 < \o(A') < \o(A)$. Of course, such a class can only exist if
the period group is discrete.

The Gromov-Witten invariant which counts rational pseudo-holomorphic 
$A$-curves with three marked points can be written as a bilinear map
$\psi_A: R \times R \longrightarrow R$. This map is (graded) symmetric,
of degree $-2\leftsc c_1\mo,A\rightsc$, and has the following properties:
\begin{normallist}
\item \label{item:psi-associativity}
$x\psi_A(y,z) - \psi_A(xy,z) + \psi_A(x,yz) - \psi_A(x,y)z = 0$,
\item \label{item:psi-unit}
$\psi_A(1,x) = 0$ for all $x$, and
\item \label{item:psi-divisor}
if $u \in R^2$ satisfies $\leftsc u,A \rightsc = 0$ then
$\psi_A(u,x) = 0$ for all $x$.
\end{normallist}
Define a product $\qpa$ on $R \otimes \field[q]/q^2$ by
\[
(x_0 + x_1q) \qpa (y_0 + y_1q) = (x_0y_0) + (x_0y_1+x_1y_0 + 
\psi_A(x_0,y_0))q.
\]
Property \ref{item:psi-associativity} is equivalent to the associativity
of $\qpa$. Property \ref{item:psi-unit} says that $1 \in R$ is a unit
for $\qpa$, and property \ref{item:psi-divisor} can be written as
$u \qpa (x_0 + x_1q) = ux_0 + (ux_1)q$.

\begin{defn} \label{def:extension} Let $(\tR,t,j)$ be a deformation of $R$. 
An extension of $\psi_A$ to $\tR$ is a bilinear map
$\tpsi_A: \tR \times \tR \longrightarrow \tR$ of degree
$-2\leftsc c_1\mo,A \rightsc$ which is (graded) symmetric and 
has the following properties:
\begin{primelist}
\item \label{item:tpsi-associativity}
$x\tpsi_A(y,z) - \tpsi_A(xy,z) + \tpsi_A(x,yz) - \tpsi_A(x,y)z = 0$;
\item \label{item:tpsi-unit}
$\tpsi_A(1,x) = \tpsi_A(t,x) = 0$ for all $x$;
\item \label{item:tpsi-divisor}
if $u \in \tR^2$ satisfies $\leftsc j(u),A \rightsc = 0$ then 
$\tpsi_A(u,x) = 0$ for all $x$;
\item \label{item:tpsi-ring-homomorphism}
$j(\tpsi_A(x,y)) = \psi_A(j(x),j(y))$.
\end{primelist}
\end{defn}

These properties can again be interpreted in terms of a suitably
defined product $\tqpa$ on $\tR \otimes \field[q]/q^2$. 
\ref{item:tpsi-ring-homomorphism} says that
\[
j \otimes \id_{\field[q]/q^2}: (\tR \otimes \field[q]/q^2,\tqpa)
\longrightarrow (R \otimes \field[q]/q^2,\qpa)
\]
is a ring homomorphism. The deformations of $R$ which admit an extension 
of $\psi_A$ form a subset $\Def_d(R,\psi_A) \subset \Def_d(R)$. 
We leave it to the reader to verify that this is actually a linear subspace.

Let $(E,\pi,i)$ be a smooth fibre bundle over $S^d$ with fibre $M$.
A fibrewise symplectic structure on $E$ is a smooth family 
$\O = (\O_b)_{b \in S^d}$ of symplectic forms on the fibres $E_b$ such that 
the cohomology class $[\O_b] \in H^2(E_b;\R)$ is locally constant 
in $b$. If in addition $i^*(\O_{b_0}) = \o$, we call
$(E,\pi,i,\O)$ a symplectic fibre bundle with fibre $\mo$.
The isomorphism classes of such bundles form a group under fibre connected 
sum, and this group is isomorphic to $\pi_{d-1}(\Aut\mo)$. Forgetting 
$\O$ corresponds to passing from $\Aut\mo$ to $\Diff(M)$.

\begin{prop} \label{th:fibrewise-quantum-product}
Let $(E,\pi,i,\O)$ be a symplectic fibre bundle over $S^d$ ($d \geq 2$) 
with fibre $\mo$, such that $i^*$ is onto. Then the deformation of $R$
determined by $H^*(E)$ admits an extension of $\psi_A$. \end{prop}

In the case when $G_{d-1} = \pi_{d-1}(\Diff(M))$, the situation can be 
summarized in the commutative diagram
\[
\xymatrix{
\pi_{d-1}(\Aut\mo) 
\ar[r] \ar[d] &
\pi_{d-1}(\Diff(M)) \ar[d]^{\alpha_{d-1}} \\
{\Def_d(R,\psi_A)} \ar@{^{(}->}[r] & {\Def_d(R)}.
}
\]

%------------------\input{5}

\section{Proof of Proposition \ref{th:fibrewise-quantum-product} 
\label{sec:five}}

We begin by recalling the definition of the Gromov-Witten invariant
$\psi_A$. Let $J$ be an $\o$-compatible almost complex structure on $M$,
$\moduli(A,J)$ the moduli space of $J$-holo\-mor\-phic maps
$u: \CP{1} \longrightarrow M$ which represent $A$, and
\[
\moduli_r(A,J) = \moduli(A,J) \times_{PSL_2(\C)} (\CP{1})^r
\]
($r = 0,1,2,\dots$) the associated moduli spaces of marked 
$J$-holomorphic curves. Our assumptions on $A$ imply that 
any $u \in \moduli(A,J)$ is simple (not multiply covered) and that 
$A$ cannot be represented by a $J$-holomorphic cusp curve. It follows
from the standard theory of pseudo-holomorphic curves (see 
\cite{mcduff-salamon} for an exposition) that for generic $J$
the spaces $\moduli_r(A,J)$ are compact oriented smooth manifolds.
Let $e_r: \moduli_r(A,J) \longrightarrow M^r$ be the $r$-fold evaluation 
map. $\psi_A$ is defined by
\[
\leftsc \psi_A(x,y)z, [M]\rightsc = \leftsc x \times y \times z,
(e_3)_*[\moduli_3(A,J)] \rightsc
\]
for $x,y,z \in H^*(M)$ and generic $J$. A cobordism argument
\cite[Theorem 3.1.3]{mcduff-salamon} shows that $\psi_A$ is independent
of the choice of $J$. 

Now let $(E,\pi,i,\O)$ be a symplectic fibre bundle over $S^d$ $(d \geq 2)$
with fibre $\mo$. The assumption that $d \geq 2$ implies that the homology
of any fibre $E_b$ can be canonically identified with the homology of 
$E_{b_0}$ and hence (using $i$) with that of $M$. Therefore it makes sense to
say that a map $u: \CP{1} \longrightarrow E_b$ represents $A \in H_2^s(M;\Z)$.
Let $\J = (J_b)_{b \in S^d}$ be a family of almost complex structures on
the fibres $E_b$ such that $J_b$ is $\O_b$-compatible for all $b$. 
Equivalently, one can think of $\J$ as an almost complex structure on the
subbundle $TE^v = \ker(T\pi) \subset TE$, with the property that 
$\O(\cdot,\J\cdot)$ is a metric on $TE^v$. Let $\moduli^p(A,\J)$ be the
space of pairs $(b,u)$ where $b \in S^d$ and $u: \CP{1} \longrightarrow E_b$
is a $J_b$-holomorphic map which represents $A$ in the sense explained above.
$\moduli^p(A,\J)$ is called the {\em parametrized moduli space} of rational 
$\J$-holomorphic $A$-curves on $E$. As before, we define
\[
\moduli_r^p(A,\J) = \moduli^p(A,\J) \times_{PSL_2(\C)} (\CP{1})^r.
\]
The basic theory of pseudo-ho\-lo\-mor\-phic curves carries 
over to the parame\-trized
situation (this is familiar from the cobordism argument which we have
mentioned above, in which one considers one-parameter families of almost
complex structures on $M$). In particular, for generic $\J$ the spaces
$\moduli_r^p(A,\J)$ are again compact smooth oriented manifolds; their
dimension is 
\[
\dim \moduli_r^p(A,\J) = \dim E + 2 \leftsc c_1\mo, A \rightsc + 2(r-3).
\]
$\moduli_r^p(A,\J)$ comes with two canonical maps: the projection
$P_r: \moduli_r^p(A,\J) \longrightarrow S^d$ and the $r$-fold evaluation
map $e_r^p: \moduli_r^p(A,\J) \longrightarrow E^r$. Note that in
general $P_r$ is not a smooth fibration. We define the {\em
parametrized Gromov-Witten invariant} $\tpsi_A: H^*(E) \times H^*(E)
\longrightarrow H^*(E)$ by
\[
\leftsc \tpsi_A(x,y)z,[E] \rightsc = \leftsc x \times y \times z,
(e_3^p)_* [\moduli_3^p(A,\J)] \rightsc
\]
for generic $\J$. Here $[E]$ is the orientation induced by the symplectic
orientation of $M$ and the standard orientation of $S^d$. An argument similar
to that for $\psi_A$ proves that $\tpsi_A$ is independent of $\J$.

The three properties of $\psi_A$ listed in section \ref{sec:four} are
special cases of properties which are known to hold for far more general
Gromov-Witten invariants. Nevertheless, we will outline proofs of them,
making use of the special features of our case to simplify the argument.

(1) Consider the cycles $Z,Z' \subset \moduli_4(A,J)$ consisting of
those $[u,z_1,z_2,z_3,z_4]$ such that $z_1 = z_2$ or $z_3 = z_4$ (for
$Z$) resp. $z_2 = z_3$ or $z_1 = z_4$ (for $Z'$).
One can prove easily, using the cross-ratio $(z_1,z_2,z_3,z_4) \longmapsto 
\frac{(z_1-z_2)(z_3-z_4)}{(z_2-z_3)(z_1-z_4)}$, that $Z$ and $Z'$ are
homologous. On the other hand
\begin{align*}
\leftsc x \times y \times z \times w, (e_4)_*[Z] \rightsc 
&= \leftsc xy \times z \times w + x \times y \times zw, 
(e_3)_* [\moduli_3(A,J)] \rightsc \\
&= \leftsc \psi_A(xy,z)w + \psi_A(x,y)zw, [M] \rightsc \\
\intertext{and similarly}
\leftsc x \times y \times z \times w, (e_4)_*[Z'] \rightsc
&= \leftsc \psi_A(x,yz)w + x\psi_A(y,z)w, [M] \rightsc.
\end{align*}

(2) Let $\phi: \moduli_3(A,J) \longrightarrow \moduli_2(A,J)$ be the
map which forgets the first marked point. There is a commutative diagram
\[
\xymatrix{
\moduli_3(A,J) \ar[rr]^{e_3} \ar[d]_{\phi}
&& M^3 \ar[d]^{p} \\
\moduli_2(A,J) \ar[rr]^{e_2}
&& M^2}
\]
where $p$ is the projection. Because $\dim \moduli_2(A,J) = \dim
\moduli_3(A,J) - 2$, $\phi_*[\moduli_3(A,J)] = 0$. This implies that
\[
\leftsc 1 \times y \times z, (e_3)_*[\moduli_3(A,J)] \rightsc
= \leftsc y \times z, p_*(e_3)_*[\moduli_3(A,J)] \rightsc = 0.
\]

(3) The forgetful map $\phi$ is a fibration with fibre $\CP{1}$.
Let $\phi_!$ be the cohomology transfer.
\begin{align*}
\leftsc \psi_A(u,y) z , [M] \rightsc &=
\leftsc e_3^*(u \times y \times z), [\moduli_3(A,J)] \rightsc \\
&= \leftsc e_3^*(u \times 1 \times 1) \phi^*(e_2^*(y \times z)), 
[\moduli_3(A,J)] \rightsc\\
& = \leftsc \phi_!(e_3^*(u \times 1 \times 1)) e_2^*(y \times z),
[\moduli_2(A,J)] \rightsc
\end{align*}
for all $u,y,z \in H^*(M)$. Now assume that $u$ is two-dimensional.
Since $e_3$ maps the homology class of the fibre of $\phi$ to 
$A \times [\point] \times [\point] \in H_2(M^3)$, one has
\[
\phi_! e_3^*(u \times 1 \times 1) = \leftsc u, A \rightsc 1 
\in H^0(\moduli_2(A,J)).
\]
Hence $\psi_A(u,y)$ vanishes if $\leftsc u, A \rightsc = 0$.

To prove Proposition \ref{th:fibrewise-quantum-product} we have to
show that $\tpsi_A$ is an extension of $\psi_A$ in the sense of
Definition \ref{def:extension}. The first three properties listed
there are analogues of the properties of $\psi_A$, and the proofs 
given above can be easily adapted to the parametrized case. The 
remaining property is
\begin{equation} \label{eq:last-property}
i^*(\tpsi_A(x,y)) = \psi_A(i^*(x),i^*(y)).
\end{equation}
As usual, let $t \in H^d(E)$ be the pullback of the fundamental class
of $S^d$. $1 \times 1 \times t \in H^d(E^3)$ is Poincar{\'e} dual
to the submanifold $C = E \times E \times E_{b_0} \subset E^3$. Therefore
\begin{align*}
\leftsc i^*(\tpsi_A(x,y))i^*(z),[M] \rightsc
& = \leftsc \tpsi_A(x,y) z, [E_{b_0}] \rightsc \\
& = \leftsc \tpsi_A(x,y) zt, [E] \rightsc \\
& = \leftsc x \times y \times zt, (e_3^p)_*[\moduli_3^p(A,\J)] \rightsc\\
& = \leftsc x \times y \times z, (e_3^p)_*[\moduli_3^p(A,\J)] \circ
[C] \rightsc 
\end{align*}
for all $x,y,z \in H^*(E)$. Here $\circ$ denotes the intersection product on
$H_*(E^3)$. Now, if we assume that
\begin{equation} \label{eq:fundamental-relation}
(e_3^p)_*[\moduli_3^p(A,\J)] \circ [C] =
(i \times i \times i)_*(e_3)_*[\moduli_3(A,J)],
\end{equation}
we obtain
\begin{multline*}
\leftsc i^*(\tpsi_A(x,y))i^*(z),[M] \rightsc =
\leftsc i^*(x) \times i^*(y) \times i^*(z), (e_3)_*[\moduli_3(A,\J)] 
\rightsc = \\
= \leftsc \psi_A(i^*(x),i^*(y)) i^*(z), [M] \rightsc.
\end{multline*}
Since $i^*(z)$ can be any element of $H^*(M)$ (we have assumed in
Proposition \ref{th:fibrewise-quantum-product} that $i^*$ is onto), this
implies \eqref{eq:last-property}. 

Equation \eqref{eq:fundamental-relation} expresses the (rather obvious) 
fact that for suitable choices of $J$ and $\J$, the unparametrized moduli
space $\moduli_3(A,J)$ can be identified with a fibre of the projection 
$P_3: \moduli_3^p(A,\J) \longrightarrow S^d$. A more precise formulation
of the argument goes as follows: choose $J$ and $\J$ in 
such a way that the moduli spaces $\moduli_3(A,J)$ and $\moduli^p_3(A,\J)$ 
are smooth and $i_*(J) = J_{b_0}$; this is possible. Then $e_3^p$ is 
transverse to $C$, and $(e_3^p)^{-1}(C) \subset \moduli_3^p(A,\J)$, which 
consists of those elements $[b,u,z_1,z_2,z_3]$ with $b = b_0$, can be 
identified with $\moduli_3(A,J)$ in such a way that the diagram
\[
\xymatrix{
\moduli_3(A,J) \ar[r]^-{e_3} \ar[d]_{=} 
& M^3 \ar[d]^{i \times i \times i} \\
(e_3^p)^{-1}(C) \ar[r]^-{e_3^p}
& E^3
}
\]
commutes. This proves \eqref{eq:fundamental-relation} and hence completes
the proof of Proposition \ref{th:fibrewise-quantum-product}.

%------------------\input{6}

\section{The main computation \label{sec:six}}

Let $R_1,R_2$ be a pair of graded $\field$-algebras, and $R = R_1 \otimes R_2$
their graded tensor product. Given a deformation $(\tR_1,t_1,j_1)$ of $R_1$
and a deformation $(\tR_2,t_2,j_2)$ of $R_2$ of the same dimension, one
can define a deformation $(\tR,t,j)$ of $R$ by
$
\tR = \tR_1 \otimes_{\field[\epsilon]/\epsilon^2} \tR_2
$
(where $\epsilon$ acts on $\tR_i$ by multiplication with $t_i$),
$t = t_1 \otimes 1 = 1 \otimes t_2$, and $j = j_1 \otimes j_2$.
This operation, which we call the exterior product of deformations, 
defines homomorphisms $\times_d:
\Def_d(R_1) \oplus \Def_d(R_2) \longrightarrow \Def_d(R)$.
We call a deformation of $R$ {\em split} if it is isomorphic 
to the exterior product of deformations of $R_i$, that is, if its isomorphism 
class lies in the image of $\times_d$.

\begin{remark} In the alternative picture of Remark \ref{th:explicit}
the exterior product is defined by assigning to a pair of 
bilinear maps $\psi_i: R_i \otimes R_i \longrightarrow R_i$ the product
\begin{multline*}
\psi(x_1 \otimes x_2,y_1 \otimes y_2) = 
(-1)^{\deg(x_2)\deg(y_1)} \big(\psi_1(x_1,y_1) \otimes x_2y_2 + \\
+ (-1)^{\deg(x_1y_1)d} x_1y_1 \otimes \psi_2(x_2,y_2) \big).
\end{multline*}
\end{remark}

The main result of this section is

\begin{prop} \label{th:symplectic-implies-split}
Let $(E,\pi,i,\O)$ be a symplectic fibre bundle over
$S^d$, for some even $d$, with fibre $(P_{mn},\eta_{mn})$. Then the
deformation of $H^*(P_{mn}) = H^*(\CP{m}) \otimes H^*(\CP{n})$
determined by $H^*(E)$ is split.
\end{prop}

We will obtain this as a consequence of a more general result. To state that
result we need to introduce one more algebraic notion. 
Let $f_i: R_i \longrightarrow R$ be the obvious inclusions.

\begin{defn} Let $(\tR,t,j)$ be a deformation of $R$. We say that
$\tR$ is {\em semi-split} with respect to $R_i$ ($i = 1$ or $2$) if there 
is a deformation $(\tR_i,t_i,j_i)$ of $R_i$ and a morphism 
$\tilde{f}_i: \tR_i \longrightarrow \tR$ which lies over $f_i$.
\end{defn}

A split deformation is clearly semi-split with respect to both
$R_1$ and $R_2$. The converse it also true: a deformation of $R$ which
is semi-split with respect to both $R_1$ and $R_2$ is split.
Later, we will use the following elementary criterion:

\begin{lemma} \label{th:invariant-subspace}
Let $(\tR,t,j)$ be a deformation of $R$. If $\tR$ has a subalgebra
$\tR_1$ with $t \in \tR_1$, $j(\tR_1) = \im(f_1)$
and $\dim (\tR_1 \cap t \tR) \leq \dim R_1$ then it is semi-split
with respect to $R_1$. \qed
\end{lemma}

Let $M'$ be a compact manifold whose cohomology ring is generated
by $H^2(M')$, and $\o'$ a symplectic form on $M'$ such that 
$\o'(H_2^s(M';\Z)) \subset \Z$. We consider the product 
$\mo = (M',\o') \times (\CP{n},\o_n)$ for some $n$. Let $A \in
H_2^s(M;\Z)$ be the homology class of a line in $\CP{n}$. $A$
satisfies the conditions of section \ref{sec:three}, and the
corresponding Gromov-Witten invariant is well-known:
\begin{equation}
\label{eq:truncated-product}
\psi_A(x \otimes u^i,y \otimes u^j) =
\begin{cases}
xy \otimes u^{i+j-n-1} & i+j \geq n+1\\
0 & \text{otherwise}
\end{cases}
\end{equation}
for all $x,y \in H^*(M')$ and $0 \leq i,j \leq n$ (as usual,
$u$ denotes the standard generator of $H^2(\CP{n})$).

\begin{prop} \label{th:symplectic-implies-semi-split} 
Let $(E,\pi,i,\O)$ be a symplectic fibre bundle with fibre $\mo$ 
over $S^d$ for some even $d$. Then the deformation of 
$R = H^*(M)$ given by $H^*(E)$ is semi-split with respect to $H^*(M')$.
\end{prop}

Setting $(M',\o') = (\CP{m},\o_m)$, it follows that a symplectic fibre
bundle over $S^d$ with fibre $(P_{mn},\eta_{mn})$ determines a deformation
of $H^*(P_{mn})$ which is semi-split with respect to $H^*(\CP{m})$. Since
the situation is symmetric with respect to $m$ and $n$, it follows that the
deformation is also semi-split with respect to $H^*(\CP{n})$, hence split.
Therefore Proposition \ref{th:symplectic-implies-semi-split} implies
Proposition \ref{th:symplectic-implies-split}. 

\proof[Proof of Proposition \ref{th:symplectic-implies-semi-split}]
Note that $i^*$ is onto because $d$ is even and $H^*(M)$ is concentrated
in even dimensions. Let $(\tR,t,j)$ be the deformation of $R$ determined
by $H^*(E)$. Proposition \ref{th:fibrewise-quantum-product} ensures the
existence of an extension $\tpsi_A$ of $\psi_A$ to $\tR$. Choose a 
basis $\xi_1,\dots,\xi_g$ of 
$H^2(M')$ and elements $\txi_1,\dots,\txi_g,\tu \in \tR^2$ such
that $j(\txi_i) = \xi_i \otimes 1$ and $j(\tu) = 1 \otimes u$. Let
$\tR_1 \subset \tR$ be the subalgebra generated by $t,\txi_1,\dots,\txi_g$. 
Clearly $j(\tR_1) = H^*(M') \otimes 1$. We want to apply Lemma
\ref{th:invariant-subspace} to this subalgebra. What remains to be shown
is that the dimension of $\tR_1 \cap t \tR$ is not greater than the 
dimension of $H^*(M')$. This is a consequence of the following

\begin{assertion*} Every $x \in \tR_1 \cap t \tR$ is of the form
$x = ty$ for some $y$ such that $j(y) \in H^*(M') \otimes 1$.
\end{assertion*}

The Assertion is proved in two steps:

{\bf Step 1:} {\em Every $x \in \tR_1$ satisfies $\tpsi_A(x,\tu^n) = 0$.}
Since $\leftsc j(\txi_i),A \rightsc = \leftsc \xi_i, A \rightsc = 0$,
$\tpsi_A(\xi_i,x) = 0$ for all $x$ by property \ref{item:tpsi-divisor} of 
$\tpsi_A$. Therefore the equation
\[
\txi_i \tpsi_A(y,\tu^n) - \tpsi_A(\txi_i y, \tu^n) +
\tpsi_A(\txi_i, y \tu^n) - \tpsi_A(\txi_i,y) \tu^n = 0,
\]
which is a special case of property \ref{item:tpsi-associativity} of 
$\tpsi_A$, simplifies to
\begin{equation} \label{eq:inductively}
\tpsi_A(\txi_i y,\tu^n) = \txi_i \tpsi_A(y,\tu^n).
\end{equation}
Note that by property \ref{item:tpsi-unit} $\tpsi_A(1,\tu^n) = 
\tpsi_A(t,\tu^n) = 0$. Arguing inductively using \eqref{eq:inductively} it 
follows that $\tpsi_A(x,\tu^n) = 0$ for all $x \in \tR_1$.

{\bf Step 2.}
{\em If $x = ty \in t\tR$ satisfies $\tpsi_A(x,\tu^n) = 0$,
then $j(y) \in H^*(M') \otimes 1$.} Since $\tpsi_A(t,w) = 0$ for
all $w$, the equation $t \tpsi_A(y,\tu^n) - \tpsi_A(x,\tu^n) + 
\tpsi_A(t,y\tu^n) - \tpsi_A(t,y)\tu^n = 0$ reduces to
\begin{equation} \label{eq:product-with-t}
\tpsi_A(x,\tu^n) = t \tpsi_A(y,\tu^n).
\end{equation}
Now write $j(y) = \sum_{i=0}^n c_i \otimes u^i$ with $c_i \in H^*(M')$, and 
choose $\tilde{c}_0, \dots, \tilde{c}_n \in \tR$ with $j(\tilde{c}_i) = c_i$. 
Using property \ref{item:tpsi-ring-homomorphism} of $\tpsi_A$ and the formula
\eqref{eq:truncated-product} for $\psi_A$ one sees that
\[
j(\tpsi_A(y,\tu^n)) = \psi_A(j(y),u^n) = \sum_{i=1}^n c_i \otimes u^{i-1}.
\]
Therefore $\tpsi_A(y,\tu^n) = \sum_{i=1}^n \tilde{c}_i \tu^{i-1} + tz$
for some $z \in \tR$. Using \eqref{eq:product-with-t} it follows that
\[
\tpsi_A(x,\tu^n) = t\Big(\sum_{i=1}^n \tilde{c}_i \tu^{i-1}\Big).
\]
If this vanishes then $t\tilde{c}_i = 0$ and hence $c_i = 0$ for all
$i>0$, which means that $j(y) \in H^*(M') \otimes 1$. \qed

%------------------\input{7}

\section{Conclusion \label{sec:seven}}

We can now prove Theorem \ref{th:main}. Set $\field = \Q$, and let 
$R = H^*(P_{mn})$. Fix an even $d$, and let $\Def^s_d(R) \subset \Def_d(R)$
be the subspace of split deformations. According to Proposition
\ref{th:symplectic-implies-split} there is a commutative diagram
\[
\xymatrix{
\pi_{d-1}(\Aut(P_{mn},\eta_{mn})) \otimes \Q
\ar[rr]^-{\beta_{d-1} \otimes \id_{\Q}} \ar[d] &&
\pi_{d-1}(\Diff(M)) \otimes \Q \ar[d]^{\alpha_{d-1} \otimes \id_{\Q}} \\
{\Def_d^s(R)} \ar@{^{(}->}[rr] && {\Def_d(R)}
}
\]
Proposition \ref{th:surjective} shows that $\alpha_{d-1} \otimes \id_{\Q}$
is onto. Therefore
\[
\rank(\mathrm{coker}\; \beta_{d-1}) \geq \dim_{\Q} \Def_d(R) -
\dim_{\Q} \Def_d^s(R).
\]
Assume first that $d>2$. Using Example \ref{ex:cpn} and the results 
of section \ref{sec:three} one can identify the subspace 
$\Def^s_d(R) \subset \Def_d(R)$ with 
$H^{2m+2-d}(\CP{m}) \oplus H^{2n+2-d}(\CP{n}) \subset 
R^{2m+2-d} \oplus R^{2n+2-d}$. Therefore 
\begin{multline*}
\rank(\coker \beta_{d-1}) \geq b_{2m+2-d}(P_{mn}) - b_{2m+2-d}(\CP{m}) + \\ 
+ b_{2n+2-d}(P_{mn}) - b_{2n+2-d}(\CP{n}). 
\end{multline*}
This is positive as long as $d \leq \max\{2m,2n\}$. For $d = 2$ one has 
$\Def_2^s(R) = 0$ (by Example \ref{ex:cpn}) and hence
$\rank(\coker \beta_1) \geq \dim_{\Q} \Def_2(R) =
\dim_{\Q} R^{2m} + \dim_{\Q} R^{2n} - 2 =
b_{2m}(P_{mn}) - b_{2m}(\CP{m}) + b_{2n}(P_{mn}) - b_{2n}(\CP{n})$.
\qed

One can ask what happens if one considers $P_{mn}$ with a weighted
symplectic form $\eta_{mn}^{(\lambda)} = \lambda(\o_m \times 1) +
1 \times \o_n$ where $\lambda>1$. Proposition 
\ref{th:symplectic-implies-split} is no longer true for symplectic fibre 
bundles with fibre $(P_{mn},\eta_{mn}^{(\lambda)})$, at least if $\lambda$
is sufficiently large. Indeed, all the examples constructed in Proposition
\ref{th:surjective}, which were of the form $E = \P(\xi)$ for some
$\xi \longrightarrow \CP{m} \times S^d$, admit a fibrewise symplectic
structure modelled on $\eta_{mn}^{(\lambda)}$ for some large $\lambda$,
while the corresponding deformations are not always split.
If $\lambda$ is integral, we can still use Proposition 
\ref{th:symplectic-implies-semi-split} to obtain some information about
symplectic fibre bundles with fibre $(P_{mn},\eta_{mn}^{(\lambda)})$. 
This leads to the following

\begin{theorem} Let $\lambda \geq 2$ be an integer. Then the homomorphism
\[
\beta_k^{(\lambda)}: 
\pi_k(\Aut(P_{mn},\eta_{mn}^{(\lambda)})) \longrightarrow \pi_k(\Diff(P_{mn}))
\]
induced by inclusion is not surjective for any odd $k$ with
$1 \leq k \leq 2m-1$. In fact
$\rank(\coker \beta_k^{(\lambda)}) \geq 
b_{2m+1-k}(P_{mn}) - b_{2m+1-k}(\CP{m}) > 0$.
\end{theorem}

We omit the proof, which is similar to that of Theorem \ref{th:main}.
For $\lambda \notin \N$ it becomes possible that the class 
$A \in H_2^s(P_{mn};\Z)$ coming from the generator of $H_2(\CP{n};\Z)$
can be represented by a rational pseudo-holomorphic cusp-curve.
However, our argument would still work if one could prove that these 
cusp-curves can be removed by perturbing the almost complex structure. 
More precisely, one needs to show that no such cusp-curves occur 
in a generic family of compatible almost complex structures depending 
on a certain number of parameters.

\begin{example}
Let $(E,\pi,i,\O)$ be a symplectic fibre bundle over $S^2$ with fibre
$(P_{21},\eta_{21}^{(\lambda)})$. Choose a family $\J = (J_b)_{b \in S^2}$
of compatible almost complex structures on its fibres. If $A$ is 
represented by a rational $J_b$-holomorphic cusp-curve for some $b$,
there must be $A_1,A_2 \in H_2^s(P_{21};\Z)$ with $A_1 + A_2 = A$, such that
$\o(A_1),\o(A_2) > 0$ and 
$\moduli^f(A_1,\J), \moduli^f(A_2,\J) \neq \emptyset$. The virtual
dimensions of these moduli spaces are
$\vdim \moduli^f(A_1,\J) = 8 + 2\leftsc c_1,A_1 \rightsc$,
$\vdim \moduli^f(A_2,\J) = 8 + 2\leftsc c_1,A_2 \rightsc$.
A pair $(A_1,A_2)$ such that both dimensions are nonnegative exists
iff there are integers $k,l$ with
$1 > \lambda k + l > 0$ and $6 \geq 3k + 2l \geq -4$.
An elementary argument shows that this is possible only if $\lambda < 3$.
Therefore it is plausible that if $\lambda \geq 3$ then the deformation 
of $H^*(P_{21})$ determined by $H^*(E)$ is semi-split with 
respect to $H^*(\CP{2})$. To make this argument rigorous, one has to 
deal with the problem of multiply-covered 
pseudo-holomorphic maps of negative Chern number, in the manner of 
\cite{li-tian96}, \cite{ruan96}, \cite{siebert96}.
\end{example}

It is interesting to compare our approach with the geometric
methods used by Gromov and others in the four-dimensional case. For
instance, consider the following result, which is a version of the main
step in the proof of Gromov's theorem on $\Aut(P_{11},\eta_{11})$.

\begin{proposition}[Gromov] \label{th:gromov-splitting}
Any symplectic fibre bundle over $S^d$ ($d \geq 2$) with fibre 
$(P_{11},\eta_{11})$ is the fibre product of two symplectic fibre bundles
with fibre $(\CP{1},\o_1)$.
\end{proposition}

Proposition \ref{th:symplectic-implies-split} is a higher-dimensional 
analogue of this in which the splitting takes place only on the
level of cohomology rings. The reason for the stronger nature of the
four-dimensional result is that the geometric behaviour of 
pseudo-holomorphic curves on a symplectic four-manifold is tightly 
controlled by the positive intersection theorem and the 
adjunction inequality \cite{mcduff93b}, both of which are used in
the proof of Proposition \ref{th:gromov-splitting}.
We end this discussion with a general conjecture suggested
by Proposition \ref{th:symplectic-implies-semi-split}.

\begin{conjecture} Let $\mo$ and $(N,\eta)$ be compact symplectic
manifolds. Let $(E,\pi,i)$ be a smooth fibre bundle over $S^d$ with
fibre $M \times N$, such that $i^*:H^*(E) \longrightarrow H^*(M \times N)$
is surjective. Assume that for all sufficiently large $\lambda$ there
is a fibrewise symplectic structure $\O^{(\lambda)}$ on $E$ such that
$i^*(\O^{(\lambda)}) = \lambda(\o \times 1) + 1 \times \eta$. Then
the deformation of $H^*(M \times N)$ determined by $H^*(E)$ is semi-split
with respect to $H^*(M)$. \end{conjecture}

%
%\bibliographystyle{/emtex/amslatex/amsplain}
%\bibliographystyle{/opt/texmf/tex/latex/packages/amslatex/classes/amsplain}
%\bibliographystyle{amsplain}
%\bibliography{../macros/standard,../macros/new}

\providecommand{\bysame}{\leavevmode\hbox to3em{\hrulefill}\thinspace}

\end{document}